\documentclass[a4paper]{article}
\usepackage{amssymb}
\usepackage{amsfonts}
\usepackage{amsmath}
\usepackage{amsthm}
\usepackage{authblk}
\usepackage[USenglish]{babel}
\usepackage{esint}
\usepackage{graphicx}
\usepackage[colorlinks=true]{hyperref}
\usepackage[latin1]{inputenc}
\newtheorem{theorem}{Theorem}[section]
\newtheorem{lemma}[theorem]{Lemma}
\newtheorem{proposition}[theorem]{Proposition}
\newtheorem{corollary}[theorem]{Corollary}
\newtheorem{remark}[theorem]{Remark}
\newtheorem{definition}[theorem]{Definition}

\numberwithin{equation}{section}
\newcommand{\N}{\mathbb N}
\newcommand{\Z}{\mathbb Z}

\newcommand{\R}{\mathbb R}

\renewcommand{\S}{\mathbb S}
\renewcommand{\P}{\mathbb P}

\newcommand{\mcal}{\mathcal}

\newcommand{\mrm}{\mathrm}

\renewcommand{\d}{\delta}
\newcommand{\D}{\Delta}
\newcommand{\e}{\varepsilon}
\newcommand{\z}{\zeta}
\renewcommand{\t}{\theta}

\newcommand{\la}{\lambda}

\newcommand{\s}{\sigma}
\newcommand{\si}{\varsigma}

\newcommand{\ph}{\varphi}

\newcommand{\wt}{\widetilde}

\newcommand{\fr}{\frac}
\newcommand{\pa}{\partial}
\newcommand{\n}{\nabla}
\newcommand{\fa}{\forall}

\newcommand{\sm}{\setminus}

\newcommand{\sub}{\subset}

\newcommand{\eq}{\equiv}

\newcommand{\x}{\times}
\renewcommand{\c}{\circ}
\newcommand{\cd}{\cdot}
\newcommand{\ds}{\dots}

\newcommand{\tx}{\text}
\newcommand{\q}{\quad}
\renewcommand{\l}{\left}
\renewcommand{\r}{\right}

\newcommand{\bthm}{\begin{theorem}}
\newcommand{\ethm}{\end{theorem}}
\newcommand{\blem}{\begin{lemma}}
\newcommand{\elem}{\end{lemma}}
\newcommand{\bprop}{\begin{proposition}}
\newcommand{\eprop}{\end{proposition}}
\newcommand{\bcor}{\begin{corollary}}
\newcommand{\ecor}{\end{corollary}}
\newcommand{\brem}{\begin{remark}}
\newcommand{\erem}{\end{remark}}
\newcommand{\bdefi}{\begin{definition}}
\newcommand{\edefi}{\end{definition}}
\newcommand{\bpf}{\begin{proof}}
\newcommand{\epf}{\end{proof}}
\newcommand{\bl}{\begin{array}{l}}
\newcommand{\bll}{\begin{array}{ll}}
\newcommand{\barr}{\begin{array}}
\newcommand{\earr}{\end{array}}
\newcommand{\bite}{\begin{itemize}}
\newcommand{\eite}{\end{itemize}}
\newcommand{\bequ}{\begin{equation}}
\newcommand{\eequ}{\end{equation}}
\newcommand{\beqa}{\begin{eqnarray}}
\newcommand{\eeqa}{\end{eqnarray}}
\newcommand{\beqy}{\begin{eqnarray*}}
\newcommand{\eeqy}{\end{eqnarray*}}
\newcommand{\bin}[2]{\left(\genfrac{}{}{0pt}{}{#1}{#2}\right)}
\newcommand{\qm}[1]{``#1''}
\allowdisplaybreaks

\begin{document}

\everymath{\displaystyle}

\title{Nonradial entire solutions for Liouville systems}
\author{Luca Battaglia\thanks{Sapienza Universit\`a di Roma, Dipartimento di Matematica, Piazzale Aldo Moro $5$, $00185$ Roma - battaglia@mat.uniroma1.it}, Francesca Gladiali\thanks{Universit\`a degli Studi di Sassari, Dipartimento Polcoming, Via Piandanna $4$, $00710$ Sassari - fgladiali@uniss.it}, Massimo Grossi\thanks{Sapienza Universit\`a di Roma, Dipartimento di Matematica, Piazzale Aldo Moro $5$, $00185$ Roma - grossi@mat.uniroma1.it}}
\date{}

\maketitle\

\begin{abstract}
\noindent We consider the following system of Liouville equations:
$$\l\{\bll-\D u_1=2e^{u_1}+\mu e^{u_2}&\tx{in }\R^2\\-\D u_2=\mu e^{u_1}+2e^{u_2}&\tx{in }\R^2\\
\int_{\R^2}e^{u_1}<+\infty,\ \int_{\R^2}e^{u_2}<+\infty.\
\earr\r.$$
We show the existence of at least $n-\left[\frac n3\right]$ global branches of nonradial solutions bifurcating from  $u_1(x)=u_2(x)=U(x)=\log\fr{64}{(2+\mu)\l(8+|x|^2\r)^2}$ at the values $\mu=-2\fr{n^2+n-2}{n^2+n+2}$ for any $n\in\N$.
\end{abstract}\

\section{Introduction}\
In this paper we consider the system
\bequ\label{sistema}
\l\{\bll-\D u_1=2e^{u_1}+\mu e^{u_2}&\tx{in }\R^2\\-\D u_2=\mu e^{u_1}+2e^{u_2}&\tx{in }\R^2\\
\int_{\R^2}e^{u_1}<+\infty,\ \int_{\R^2}e^{u_2}<+\infty.\
\earr\r.
\eequ
When $\mu=-1$ we recover the well known Toda system, whose solutions were completely classified in \cite{jw}. This type of systems received a growing interest
in recent years. We mention the paper \cite{pt} for a detailed description of the physical applications and the extensive list of references. We point out that the known results for \eqref{sistema} when $\mu\neq -1$ concern the existence (or nonexistence) of {\em radial} solutions. An important case was considered in \cite{CK} and \cite{csw} where the authors study the cooperative system
$$\l\{\bll-\D u_1=a_{11}e^{u_1}+a_{12}e^{u_2}&\tx{in }\R^2
\\-\D u_2=a_{21} e^{u_1}+a_{22}e^{u_2}&\tx{in }\R^2\\
\int_{\R^2}e^{u_1}<+\infty,\ \int_{\R^2}e^{u_2}<+\infty.\
\earr\r.$$
and under the assumption that the matrix $A=(a_{ij})_{i,j=1,2}$ is symmetric, irreducible and with nonnegative entries they prove that all solutions are radial (see also \cite{cz}).\\
If some of the entries $a_{ij}$ are negative the previous result is not true anymore. Indeed if $\mu=-1$ the classification result of Jost and Wang shows the existence of  {\em nonradial} solutions. It is the only paper where the existence of such solutions appears. Indeed existence results of {\em nonradial} solutions are a quite difficult task because the system has not a variational structure. Furthermore the lack of compactness due to the scaling and translation  invariance causes additional difficulties. 
In this paper we investigate the nonradial case by using the {\em  bifurcation theory}.\\
First note that for any $\mu>-2$ problem \eqref{sistema} always admits the branch of trivial solutions $(U_{\mu,\d,y},U_{\mu,\d,y})$ where $U_{\mu,\d,y}$ is the solution of the Liouville equation
$$-\D u=(2+\mu)e^u\quad\tx{in }\R^2,$$
namely
\bequ
U_{\mu,\d,y}:=\log\fr{64\d}{(2+\mu)\l(8\d+|\cd-y|^2\r)^2}.
\eequ
For sake of simplicity we consider
$$U_\mu:=U_{\mu,1,0}$$
and we look for solutions bifurcating from $(U_\mu,U_\mu)$. Note that, using the pair $(U_{\mu_1},U_{\mu_2})$, our technique allows to consider more general systems like
$$\l\{\bll-\D u_1=2e^{u_1}+\alpha e^{u_2}&\tx{in }\R^2\\-\D u_2=\beta e^{u_1}+2e^{u_2}&\tx{in }\R^2\\
\int_{\R^2}e^{u_1}<+\infty,\ \int_{\R^2}e^{u_2}<+\infty\
\earr\r.$$
with $\alpha=\frac{(2+\mu_2)(\mu_1^2+2\mu_1+2)}{(2+\mu_1)}$ and $\beta=\frac{(2+\mu_1)(\mu_2^2+2\mu_2+2)}{(2+\mu_2)}$.\\
We let to the interested reader the extension of our results to this case.\\
The main tool to find solutions bifurcating from $(U_\mu,U_\mu)$ will be the classical Crandall-Rabinowitz Theorem (see Theorem \ref{cr}). It is well known that this theorem has had enormous applications in many problems in mathematical analysis. However, in our case, its applicability is not at all trivial. This is mainly due to the following reasons,
\begin{itemize}
\item The geometrical invariances of the problem induce a lack of compactness
\item It is not trivial to prove that the solutions possibly founded in this way are {\em nonradial}
\item The kernel of the natural associated functional is not one dimensional
\end{itemize}
The first problem will be overcome looking for solutions with suitable invariance. This trick allows to kill the solutions of the linearized system generated by the scaling and translations invariance of the problem and to consider only the degeneracy due to the existence of other nonequivalent solutions. Further since the setting of the problem in $\R^2$ does not fit in the usual Sobolev spaces we prefer to define the operator in the unit sphere $\S^2$. 
Actually we will consider the system
\bequ\label{b1}
\l\{\bll-\D_{\S^2}\phi_1=2\l(e^\fr{\phi_1+\phi_2}2+e^\fr{\phi_1-\phi_2}2-2\r)&\tx{in }\S^2\\-\D_{\S^2}\phi_2=2\fr{2-\mu}{2+\mu}\l(e^\fr{\phi_1+\phi_2}2-e^\fr{\phi_1-\phi_2}2\r)&\tx{in }\S^2.\earr\r.
\eequ
Here we can use more standard Hilbert spaces and we also avoid problems with the behavior of the solutions at infinity.\\
Moreover using the standard stereographic projection we have that solutions in the sphere $\S^2$ eventually satisfy system  \eqref{sistema} .\\ 
The second obstruction is more delicate. In \cite{ggw} has been proved that the linearization of \eqref{sistema} at $(U_\mu,U_\mu)$, namely
\bequ
\l\{\bll-\D v_1=2e^{U_\mu}v_1+\mu e^{U_\mu}v_2&\tx{in }\R^2\\
-\D v_2=\mu e^{U_\mu}v_1+2 e^{U_\mu}v_2&\tx{in }\R^2
\earr\r.
\eequ
 admits the sequence of eigenvalues
\bequ\label{mun}
\mu_n=-2\fr{n^2+n-2}{n^2+n+2}
\eequ
and a corresponding $(2n+4)$-dimensional kernel. This kernel reduces to a $1$-dimensional one if we restrict to the space of the radial function (once we overcome the degeneracy due to the scaling invariance) and this was one of the ideas in \cite{ggw} to get branches of radial solutions bifurcating by $(U_\mu,U_\mu)$. In order to get a one-dimensional {\em nonradial} kernel the argument is more subtle. In order to give an idea let us argue as follows: let us consider
\bequ\label{i1}
\l\{\bll u_1=U_\mu+\frac{\Phi_1+\Phi_2}2\\
u_2=U_\mu+\frac{\Phi_1-\Phi_2}2
\earr\r.
\eequ
so that \eqref{sistema} becomes
\bequ\label{i2}
\l\{\bll -\D\Phi_1=(2+\mu)\l(e^{\fr{\Phi_1+\Phi_2}2+U_\mu}+
e^{\fr{\Phi_1-\Phi_2}2+U_\mu}-2e^{U_\mu}\r)&\tx{in }\R^2\\
 -\D\Phi_2=(2-\mu)\l(e^{\fr{\Phi_1+\Phi_2}2+U_\mu}-
e^{\fr{\Phi_1-\Phi_2}2+U_\mu}\r)&\tx{in }\R^2
\earr\r.
\eequ
Using polar coordinates let us consider the functions invariant for the following symmetries
$$
\mcal X=\{(\Phi_1,\Phi_2)\hbox{ such that }\Phi_1(\t,r)=\Phi_1(\t+\pi,r),\ \Phi_2(\t,r)=-\Phi_2(\t+\pi,r)\}
$$
A crucial remark is that the system \eqref{i2} is invariant under the action in $\mcal X$ and so $\mcal X$ is a natural constraint. This invariance with respect to {\em odd} mappings does not take place directly in \eqref{sistema} but it appears as we introduce the new function $\Phi_1$ and $\Phi_2$ in \eqref{i1}. Of course radial functions cannot belong to $\mcal X$ and so only nonradial functions are allowed. However the space $\mcal X$ will be not our finale choice because in general the linearization of \eqref{i2} is not $1$-dimensional. We need to consider a more complicate space involving more symmetries. This will be discussed in Section \ref{s3}.\\
Let us stress that the choice of the symmetries and antisymmetries allows from one side, to get through with the degeneracy induced by dilations and translations and from the other side it reduces the dimension of the kernel. It also allows to find different nonequivalent solutions to problem \eqref{sistema}.  
An important role in our results is played  by  the associated Legendre polynomials  $P_n^m(z)$  defined as
  \begin{equation}\label{i3}
  P_n^m(z)=\frac{(-1)^m}{2^nn!}(1-z^2)^\frac m2\frac{d^{n+m}}{dz^{n+m}}(z^2-1)^n
  \end{equation}
for $m,n\in \N$, $m\leq n$. 
Now we are in position to state our main result,
\bthm\label{teo}${}$\\
Let $n\in\N$ and $\mu_n$ be defined in \eqref{mun}. The points $(\mu_n,U_{\mu_n},U_{\mu_n})$ are non-radial multiple bifurcation points for the curve of solutions $(\mu,U_\mu,U_\mu)$. In particular, for any $m\in\N$ satisfying $\frac n3<m\le n$ there exists $\e_0>0$ and, for $\e\in(-\e_0,\e_0)$, a $C^1$-family of solutions $(\mu_n(\e),u_{1,\e},u_{2,\e})$ of \eqref{sistema}, which satisfy
\bequ\label{i20}
\l\{\bl\mu_n(0)=\mu_n\\u_{1,\e}(x)=U_{\mu_n}(x)+\e P_n^m\l(\fr{8-|x|^2}{8+|x|^2}\r)\cos(m\t)+\e(Z_{1,\e}(x)+Z_{2,\e}(x))\\
u_{2,\e}(x)=U_{\mu_n}(x)-\e P_n^m\l(\fr{8-|x|^2}{8+|x|^2}\r)\cos(m\t)+\e(Z_{1,\e}(x)-Z_{2,\e}(x))\\Z_{1,0}=Z_{2,0}=0\earr\r.,
\eequ
where $Z_{1,\e},Z_{2,\e}$ are uniformly bounded functions in $L^\infty\l(\R^2\r)$ satisfying
\bequ\label{z}
\barr{lllll}
Z_{1,\e}\l(\fr{8}re^{i\t}\r)=Z_{1,\e}\l(re^{i\t}\r)&&Z_{1,\e}\l(re^{i\l(-\t+\fr{\pi}m\r)}\r)=Z_{1,\e}\l(re^{i\t}\r)&&Z_{1,\e}\l(re^{i\l(\t+\fr{\pi}m\r)}\r)=Z_{1,\e}\l(re^{i\t}\r)\\
Z_{2,\e}\l(\fr{8}re^{i\t}\r)=(-1)^{n+m}Z_{2,\e}\l(re^{i\t}\r)&&Z_{2,\e}\l(re^{i\l(-\t+\fr{\pi}m\r)}\r)=-Z_{2,\e}\l(re^{i\t}\r)&&Z_{2,\e}\l(re^{i\l(\t+\fr{\pi}m\r)}\r)=-Z_{2,\e}\l(re^{i\t}\r).\\
\earr
\eequ
Moreover the bifurcation is global and the Rabinowitz alternative holds.\\
Finally we have the following mass-quantization property
\bequ\label{i4}
\int_{\R^2}e^{u_{1,\e}}=\int_{\R^2}e^{u_{2,\e}}=\frac{8\pi}{2+\mu_n(\e)}.
\eequ
\ethm

\begin{remark}
From the symmetry properties \eqref{z} satisfied by $Z_{1,\e},Z_{2,\e}$ we deduce that $u_{1,\e},u_{2,\e}$ verify:
$$\barr{lll}
u_{1,\e}\l(\fr{8}re^{i\t}\r)=\log\fr{r^4}{64}+\l\{\bll u_{1,\e}\l(re^{i\t}\r)&\tx{if }n+m\tx{ is even}\\u_{2,\e}\l(re^{i\t}\r)&\tx{if }n+m\tx{ is odd}\earr\r.&&u_{2,\e}\l(\fr{8}re^{i\t}\r)=\log\fr{r^4}{64}+\l\{\bll u_{2,\e}\l(re^{i\t}\r)&\tx{if }n+m\tx{ is even}\\u_{1,\e}\l(re^{i\t}\r)&\tx{if }n+m\tx{ is odd}\earr\r.\\
u_{1,\e}\l(re^{i\l(\t+\fr{\pi}m\r)}\r)=u_{2,\e}\l(re^{i\t}\r)&&u_{2,\e}\l(re^{i\l(\t+\fr{\pi}m\r)}\r)=u_{1,\e}\l(re^{i\t}\r)\\
u_{1,\e}\l(re^{i\l(-\t+\fr{\pi}m\r)}\r)=u_{2,\e}\l(re^{i\t}\r)&&u_{2,\e}\l(re^{i\l(-\t+\fr{\pi}m\r)}\r)=u_{1,\e}\l(re^{i\t}\r)
\earr$$
In particular, $u_{i,\e}\l(re^{i\l(\t+\fr{2\pi}m\r)}\r)=u_{i,\e}\l(re^{i\t}\r)=u_{i,\e}\l(re^{-i\t}\r)$ for both $i=1,2$.
\end{remark}

By \eqref{i20} and the boundedness of $Z_{1,\e},Z_{2,\e}$ we get the following
\bcor
We have that
\bequ\label{i21}
u_{1,\e}(x),u_{2,\e}(x)\sim-4\log|x|\quad\hbox{at }\infty.
\eequ
\ecor
\bcor\label{cor2}
There exists at least $n+1-\left[\frac n3\right]$ nonequivalent solutions bifurcating from $(\mu_n,U_{\mu_n},U_{\mu_n})$.
\ecor
We point out that if $n=2$ (the classical Toda system) from Theorem \ref{teo} it is possible to recover, at least locally near $(\mu_2,U_{\mu_2},U_{\mu_2})$ the classification result by Jost and Wang (see Remark \ref{e1}).\\
The global bifurcation result in Theorem \ref{teo} is new also in the case of radial solutions. In the paper \cite{ggw} only a local radial bifurcation result was proved for system \eqref{sistema}. Further in the radial setting we can separate the branches getting the following result:
\bcor\label{cor3}
The radial branches bifurcating from $(\mu_n,U_{\mu_n},U_{\mu_n})$ do not intersect.
\ecor

Finally  in the appendix we study the behavior of the parameter $\mu_n(\e)$ in a neighborhood of $\mu_{n}$. Recall that if $n=2$ by the results in \cite{jw} and \cite{wzz} we get that the set of solutions is a nondegenerate $8$-dimensional manifold. Hence in this case $\mu_2(\e)=\mu_{2}$ for any $\e$ small enough. On the other hand we will show that for many integers $m\ge1$ and $n\ge3$ this does not occur (see table at the end of appendix). Then the bifurcation diagram in a neighborhood of $\mu_{n}$ looks like very different from that of the classical Toda system. \\
The paper is organized as follows. In Section \ref{s1} we recall some preliminaries and introduce the functional setting in the sphere $\S^2$. Section \ref{s2} is devoted to the linearization of our functional and in Section \ref{s3} we introduce the symmetries described in \eqref{z}. In Section \ref{s4} we prove Theorem \ref{teo} and in the Appendix we study the behavior of the parameter $\mu_n(\e)$.

\section{Preliminaries and functional setting}\label{s1}
As mentioned in the introduction we will prove our results in the unit sphere $\S^2$. \\
A useful coordinate system on $\S^2$ will be the following hybrid between spherical and Cartesian,
$$\S^2=\l\{\l(\sqrt{1-z^2}\cos\t,\sqrt{1-z^2}\sin\t,z\r):\,z\in[-1,1],\,\t\in[0,2\pi]\r\}.$$
This coordinate system will allow to write the eigenfunctions of  $-\D_{\S^2}$  in an easier way. So we have this classical result.
\blem\label{linscal}${}$\\
The eigenvalues of the scalar linear problem
\bequ\label{scal}
\bll-\D_{\S^2} w=\la w&\tx{in }\S^2\earr
\eequ
are given by $\la_n=n(n+1)$.\\
Each eigenvalue $\la_n$ has multiplicity $2n+1$ and a basis for its eigenspace is given by
\bequ\label{base}
\{P_n^0(z),P_n^m(z)\cos(m\t),P_n^m(z)\sin(m\t)\}_{m=1,\ds,n}.
\eequ
(see \eqref{i3} for the definition of $P_n^m(z)$).
\elem
Now let us recall the celebrated result of Crandall and Rabinowitz. It will be the main tool of the paper.
\bthm{(\cite{cr}, Theorem $1.7$)}\label{cr}\\
Let $X,Y$ be Hilbert spaces, $V\sub X$ a neighborhood of $0$ and $F:(-2,2)\x V\to Y$ a map with the following properties:
\bite
\item $F(t,0)=0$ for any $t$.
\item The partial derivatives $\pa_tF$, $\pa_xF$ and $\pa^2_{t,x}F$ exist and are continuous.
\item $\ker(\pa_xF(0,0))=\mrm{span}(w_0)$ and $Z:=\mrm R(\pa_xF(0,0))^\perp$ have dimension $1$ 
\item $\pa^2_{t,x}F(0,0)w_0\notin\mrm R(\pa_xF(0,0))$.
\eite
Then, there exists $\e_0>0$, a neighborhood $U\sub(-2,2)\x X$ of $0$ and continuously differentiable maps $\eta:(-\e_0,\e_0)\to\R$ and $z:(-\e_0,\e_0)\to Z$ such that
$$\l\{\bl\eta(0)=0\\z(0)=0\\F^{-1}(0)\cap U\sm((-2,2)\x\{0\})=\{(\eta(\e),\e w_0+\e z(\e));\,\e\in(-\e_0,\e_0)\}\earr\r.$$
\ethm
Now let us introduce our functional to apply Theorem \eqref{cr}.
Let $-\D_{\S^2}$ the Laplace-Beltrami operator on $\S^2$ and $\mcal T:(-2,2)\x W^{2,2}\l(\S^2\r)\x W^{2,2}\l(\S^2\r)\to L^2\l(\S^2\r)\x L^2\l(\S^2\r)$ be defined by
\bequ\label{t}\mcal T:(\mu,\phi_1,\phi_2)\to\l(\barr{c}\D_{\S^2}\phi_1+2\l(e^\fr{\phi_1+\phi_2}2+e^\fr{\phi_1-\phi_2}2-2\r)\\
\D_{\S^2}\phi_2+2\fr{2-\mu}{2+\mu}\l(e^\fr{\phi_1+\phi_2}2-e^\fr{\phi_1-\phi_2}2\r)\earr\r).\eequ
It is not difficult to see that $\mcal T$ is  $C^2$.
Next proposition shows that zeros of   $\mcal T$ provide solutions to 
\eqref{sistema}. 
\bprop\label{sfera}${}$\\
Let $\Pi:\S^2\sm\{0,0,-1\}\to\R^2$ be 
\bequ\label{pi}
\Pi:\l(\sqrt{1-z^2}\cos\t,\sqrt{1-z^2}\sin\t,z\r)\to\l(\sqrt{8\fr{1-z}{1+z}}\cos\t,\sqrt{8\fr{1-z}{1+z}}\sin\t\r).
\eequ
a stereographic projection.
Then, for any $(\mu,\phi_1,\phi_2)$ solving $\mcal T(\mu,\phi_1,\phi_2)=\l(\barr{c}0\\0\earr\r)$ the triplet
\bequ\label{u}
(\mu,u_1,u_2)=\l(\mu,\fr{\phi_1\c\Pi^{-1}+\phi_2\c\Pi^{-1}}2+U_\mu,\fr{\phi_1\c\Pi^{-1}-\phi_2\c\Pi^{-1}}2+U_\mu\r)
\eequ
solves \eqref{sistema}.
\eprop\

\bpf${}$\\
If $(\mu,\phi_1,\phi_2)$ is a zero of $\mcal T$, then by definition it solves \eqref{b1}.\\
Let us now consider $(\Phi_1,\Phi_2)=\l(\phi_1\c\Pi^{-1},\phi_2\c\Pi^{-1}\r)$. $\Pi$ is a stereographic projection obtained projecting from the south pole $(0,0,-1)$ to the plane $\{z=\sqrt8-1\}$, therefore it is conformal; its inverse is given by
$$\Pi^{-1}:(x_1,x_2)\to\l(\fr{2\sqrt8x_1}{8+|x|^2},\fr{2\sqrt8x_2}{8+|x|^2},\fr{8-|x|^2}{8+|x|^2}\r),$$
hence the conformal factor is $\rho(x)=\fr{32}{\l(8+|x|^2\r)^2}=\fr{2+\mu}2e^{U_\mu(x)}$. Therefore, $\Phi_1$ satisfies in $\R^2$
\beqa
\label{phi1}-\D\Phi_1&=&\rho\l((-\D_{\S^2}\phi_1)\c\Pi^{-1}\r)\\
\nonumber&=&\fr{2+\mu}2e^{U_\mu}2\l(e^\fr{\phi_1\c\Pi^{-1}+
\phi_2\c\Pi^{-1}}2+e^\fr{\phi_1\c\Pi^{-1}-\phi_2\c\Pi^{-1}}2-2\r)\\
\nonumber&=&(2+\mu)\l(e^{\fr{\Phi_1+\Phi_2}2+U_\mu}+
e^{\fr{\Phi_1-\Phi_2}2+U_\mu}-2e^{U_\mu}\r)
\eeqa
and similarly $\Phi_2$ satisfies in $\R^2$
\bequ\label{phi2}
-\D\Phi_2=(2-\mu)\l(e^{\fr{\Phi_1+\Phi_2}2+U_\mu}-e^{\fr{\Phi_1-\Phi_2}2+U_\mu}\r).
\eequ
Now, take $(u_1,u_2):=\l(\fr{\Phi_1+\Phi_2}2+U_\mu,\fr{\Phi_1-\Phi_2}2+U_\mu\r)$. By \eqref{phi1} and \eqref{phi2}, $u_1$ solves
\beqy
-\D u_1=2e^{u_1}+\mu e^{u_2},
\eeqy
and we similarly get $-\D u_2=2e^{u_2}+\mu e^{u_1}$.\\
Finally, one can see that the integrability condition in \eqref{sistema} is satisfied by re-writing the integral on the sphere via the map $\Pi^{-1}$:
$$\int_{\R^2}\l(e^{u_1}+e^{u_2}\r)=\int_{\R^2}e^{U_\mu}\l(e^\fr{\Phi_1+\Phi_2}2+e^\fr{\Phi_1-\Phi_2}2\r)=\fr{2}{2+\mu}\int_{\S^2}\l(e^\fr{\phi_1+\phi_2}2+e^\fr{\phi_1-\phi_2}2\r)\le Ce^\fr{\|\phi_1\|_{L^\infty\l(\S^2\r)}+\|\phi_2\|_{L^\infty\l(\S^2\r)}}2<+\infty.$$

\epf

\section{The linearized operator}\label{s2}
In this section we study the derivative of the operator $\mcal T$. Next lemma characterizes its kernel.
\blem\label{ker}${}$\\
The linearized operator of $\mcal T$ in $(\mu,0,0)$ is given by
\bequ\label{lin}
\pa_\phi\mcal T(\mu,0,0):\l(\barr{c}w_1\\w_2\earr\r)\to\l(\barr{c}\D_{\S^2} w_1+2w_1\\\D_{\S^2} w_2+2\fr{2-\mu}{2+\mu}w_2\earr\r).
\eequ
Let $\mu_n$ be as defined in \eqref{mun}.\\ 
If $\mu\ne\mu_n$ for any $n\in\N$, then its kernel is $3$-dimensional and a basis is given by
$$\l(\barr{c}w_1(\t,z)\\w_2(\t,z)\earr\r)=span\l\{\l(\barr{c}z\\0\earr\r),\l(\barr{c}\sqrt{1-z^2}\cos\t\\0\earr\r),\l(\barr{c}\sqrt{1-z^2}\sin\t\\0\earr\r)\r\}.$$
If $\mu=\mu_n$ for some $n\in\N$, then its kernel is $2n+4$-dimensional and a basis is given by
\beqy
\l(\barr{c}w_1(\t,z)\\w_2(\t,z)\earr\r)&=&span\l\{\l(\barr{c}z\\0\earr\r),\l(\barr{c}\sqrt{1-z^2}\cos\t\\0\earr\r),\l(\barr{c}\sqrt{1-z^2}\sin\t\\0\earr\r),\l(\barr{c}0\\P_n^0(z)\earr\r),\r.\\
&&\l.\l(\barr{c}0\\P_n^m(z)\cos(m\t)\earr\r),\l(\barr{c}0\\P_n^m(z)\sin(m\t)\earr\r)\r\}_{m=1,\ds,n}
\eeqy
where $P_n^m$ are the associated Legendre polynomials as defined in \eqref{i3}.
\elem\

\bpf${}$\\
It is immediate to see, by the definition \eqref{t} of $\mcal T$, that its derivative in $\phi_1=\phi_2=0$ is given by \eqref{lin}.\\
So the kernel of $\pa_\phi\mcal T(\mu,0,0)$ is given by the solutions of
$$\l\{\bll-\D_{\S^2} w_1=2w_1&\tx{in }\S^2\\-\D_{\S^2} w_2=2\fr{2-\mu}{2+\mu}w_2&\tx{in }\S^2\earr\r.,$$
namely a system of two decoupled scalar linear problem, each of which can be solved using Lemma \ref{linscal}.\\
The first equation corresponds to \eqref{scal} with $\la=\la_1=2$ and it has a three-dimensional kernel given by
\beqy
w_1(\t,z)&=&span\l\{z,\sqrt{1-z^2}\cos\t,\sqrt{1-z^2}\sin\t\r\}.
\eeqy
The second equation has non-trivial solutions if and only if $2\fr{2-\mu}{2+\mu}=\la_n=n(n+1)$ for some $n$, namely $\mu=\mu_n$, in which case their expression is given by \eqref{base} and the dimension of the kernel is $2n+1$.\\
By combining the solutions of the two equations we get a $3$-dimensional kernel if $\mu\ne\mu_n$ and a $2n+4$-dimensional kernel if $\mu=\mu_n$ for some $n$ and the proof is concluded.
\epf\
Next lemma will be crucial to verify the transversality condition in Theorem \ref{cr}.
\blem\label{trasv}${}$\\
Assume $(w_1,w_2)\in\ker\pa_\phi\mcal T(\mu_n,0,0)$ and $\pa^2_{\mu,\phi}\mcal T(\mu_n,0,0)(w_1,w_2)\in\mrm R\l(\pa_\phi\mcal T(\mu_n,0,0)\r)$.\\
Then,
\bequ\label{span}
\l(\barr{c}w_1(\t,z)\\w_2(\t,z)\earr\r)\in\mrm{span}\l\{\l(\barr{c}z\\0\earr\r),\l(\barr{c}\sqrt{1-z^2}\cos\t\\0\earr\r),\l(\barr{c}\sqrt{1-z^2}\sin\t\\0\earr\r)\r\}
\eequ
\elem\

\bpf${}$\\
A straightforward computation shows that
\bequ\label{dmix}
\pa^2_{\mu,\phi}\mcal T(\mu_n,0,0):\l(\barr{c}w_1\\w_2\earr\r)\to\l(\barr{c}0\\-\fr{8}{(2+\mu_n)^2}w_2\earr\r).
\eequ
By the very well-known properties of $-\D_{\S^2}$ we have that
\beqy
\mrm R(\pa_\phi\mcal T(\mu_n,0,0))&:=&\l\{(\psi_1,\psi_2)\in L^2\l(\S^2\r)\x L^2\l(\S^2\r):\,\int_{\S^2}(\psi_1w_1+\psi_2w_2)=0,\,\fa\,(w_1,w_2)\in\ker(\pa_\phi\mcal T(\mu_n,0,0))\r\}.
\eeqy

If we now take $(w_1,w_2)\in\ker(\pa_\phi\mcal T(\mu_n,0,0))$ and we verify the orthogonality condition on $\pa^2_{\mu,\phi}\mcal T(\mu_n,0,0)\l(\barr{c}w_1\\w_2\earr\r)$, from \eqref{dmix} we have $(\psi_1,\psi_2)=\l(0,-\fr{8}{(2+\mu_n)^2}w_2\r)$. Hence we get
$$\int_{\S^2}(\psi_1w_1+\psi_2w_2)=-\fr{8}{(2+\mu_n)^2}\int_{\S^2}(w_2)^2=0.$$
Such a condition can only be satisfied if $w_2=0$. This leads to consider  $\pa^2_{\mu,\phi}\mcal T(\mu_n,0,0)\l(\barr{c}w_1\\0\earr\r)=\l(\barr{c}0\\0\earr\r)$ which of course belong to the range of
$\pa_\phi\mcal T(\mu_n,0,0).$ Since $\l(\barr{c}w_1\\0\earr\r)\in\ker\pa_\phi\mcal T(\mu_n,0,0)$, thanks to the characterization of the kernel given by Lemma \ref{ker} we get the claim.
\epf\

\section{The role of symmetries}\label{s3}
As pointed out in the previous section the kernel of $\pa_\phi\mcal T(\mu_n,0,0)$ is never one dimensional. So, in order to apply Theorem \ref{cr} we need to select a subspace of $W^{2,2}\l(\S^2\r)\x W^{2,2}\l(\S^2\r)$ where this property holds. It will consist of functions verifying some  symmetry and antisymmetry assumptions. So let us introduce the following  isometries,
$$\s:z\to-z\q\q\q\q\q\q\rho_m:\t\to\t+\fr{\pi}m\q\q\q\q\q\q\tau_m:\t\to-\t+\fr{\pi}m$$
where $m\geq 1$ is an integer and the following spaces,
\bequ\label{xnm}
\barr{llll}
\mcal X_{n,m}:=\{(\phi_1,\phi_2)\in W^{2,2}\l(\S^2\r)\x W^{2,2}\l(\S^2\r):&\phi_1\c\s=\phi_1,&\phi_1\c\rho_m=\phi_1,&\phi_1\c\tau_m=\phi_1,\\
&\phi_2\c\s=(-1)^{n+m}\phi_2,&\phi_2\c\rho_m=-\phi_2,&\phi_2\c\tau_m=-\phi_2\}\\
\earr
\eequ

\bequ\label{ynm}
\barr{llll}
\mcal Y_{n,m}:=\{(\psi_1,\psi_2)\in L^2\l(\S^2\r)\x L^2\l(\S^2\r):&\psi_1\c\s=\psi_1,&\psi_1\c\rho_m=\psi_1,&\psi_1\c\tau_m=\psi_1,\\
&\psi_2\c\s=(-1)^{n+m}\psi_2,&\psi_2\c\rho_m=-\psi_2,&\psi_2\c\tau_m=-\psi_2\}\\
\earr
\eequ
Note that in general $\t+\fr{\pi}m\not\in[0,2\pi]$. If this occurs we mean that in \eqref{xnm} $\phi_1\l(\t+\fr{\pi}m,z\r)=\phi_1\l(\t+\fr{\pi}m-2\pi,z\r)$ and so the definition is well posed. Same remark applies to $\tau_m$.
\brem\label{simmetrie}${}$\\
It may be interesting to consider the effect produced by the symmetries $\s,\rho_m,\tau_m$ on the plane $\R^2$ via the stereographic projection $\Pi$ defined in \eqref{pi}.\\
By parameterizing $\R^2$ using standard polar coordinates $(r,\t)$ we see that the angular variable $\t$ on $\R^2$ corresponds precisely to the variable $\t$ on $\S^2$ (hence we allow ourself a slight abuse of notation). Therefore composing $\Pi$ and $\rho_m$ will still give a rotation whereas composing $\Pi$ and $\tau_m$ will still give an angular reflection:
$$\wt\rho_m:=\Pi\c\rho_m\c\Pi^{-1}:\t\to\t+\fr{\pi}m\q\q\q\q\q\q\wt\tau_m:=\Pi\c\tau_m\c\Pi^{-1}:\t\to-\t+\fr{\pi}m.$$
On the other hand, $\wt\s:=\Pi\c\s\c\Pi^{-1}$ is no longer an isometry on the plane but rather a radial inversion with respect to the sphere of radius $\sqrt8$:
$$\wt\s:r\to\fr{8}r.$$
If we now consider how the space $\mcal X_{n,m}$ is transformed by such symmetries, namely which symmetries will inherit the solutions on the plane from the corresponding solutions on the sphere, then we get that $(\Phi_1,\Phi_2):=\l(\phi_1\c\Pi^{-1},\phi_2\c\Pi^{-1}\r)$ satisfies
$$\barr{lllll}
\Phi_1\l(\fr{8}re^{i\t}\r)=\Phi_1\l(re^{i\t}\r)&&\Phi_1\l(re^{i\l(-\t+\fr{\pi}m\r)}\r)=\Phi_1\l(re^{i\t}\r)&&\Phi_1\l(re^{i\l(\t+\fr{\pi}m\r)}\r)=\Phi_1\l(re^{i\t}\r)\\
\Phi_2\l(\fr{8}re^{i\t}\r)=(-1)^{n+m}\Phi_2\l(re^{i\t}\r)&&\Phi_2\l(re^{i\l(-\t+\fr{\pi}m\r)}\r)=-\Phi_2\l(re^{i\t}\r)&&\Phi_2\l(re^{i\l(\t+\fr{\pi}m\r)}\r)=-\Phi_2\l(re^{i\t}\r).
\earr$$
\erem\
Next lemma is one of the crucial steps that enable us to apply Theorem \ref{cr}.
\blem\label{simm}${}$\\
The restriction $\mcal T_{n,m}:=\mcal T|_{(-2,2)\x\mcal X_{n,m}}$ maps its domain into $\mcal Y_{n,m}$.\\
Moreover, if $m>\frac{n}3$, the kernel of $\pa_\phi\mcal T_{n,m}(\mu_n,0,0)$ is $1$-dimensional and it is generated by
$$\l(\barr{c}0\\P_n^m(z)\cos(m\t)\earr\r).$$
\elem\
\bpf${}$
The Laplace-Beltrami operator on $\S^2$ in the coordinates $(\t,z)$ is given by
\bequ\label{lp}
\D_{\S^2}\phi_1(\t,z)=(1-z^2)\frac{\partial^2\phi_1}{\partial z^2}-2z\frac{\partial \phi_1}{\partial z}+\frac1{1-z^2}\frac{\partial^2\phi_1}{\partial\t^2}
\eequ
By the previous formula one can easily verify that  $\mcal T=\l(\barr{c}\mcal T_1\\\mcal T_2\earr\r)$ is invariant under the reflection $\s$, i.e. $\mcal T(\mu,\phi_1\c\s,\phi_2\c\s)=\l(\barr{c}(\mcal T_1(\mu,\phi_1,\phi_2))\c\s\\(\mcal T_2(\mu,\phi_1,\phi_2))\c\s\earr\r)$ and  $\mcal T(\mu,\phi_1\c\s,-\phi_2\c\s)=\l(\barr{c}(\mcal T_1(\mu,\phi_1,\phi_2))\c\s\\-(\mcal T_2(\mu,\phi_1,\phi_2))\c\s\earr\r)$.\\
This proves that $\mcal T(\mu,\cd)$ maps the space
$$\mcal X_\pm:=\l\{(\phi_1,\phi_2)\in W^{2,2}\l(\S^2\r)\x W^{2,2}\l(\S^2\r):\,\phi_1\c\s=\phi_1,\,\phi_2\c\s=\pm\phi_2\r\}$$
into
$$\mcal Y_\pm:=\l\{(\psi_1,\psi_2)\in L^2\l(\S^2\r)\x L^2\l(\S^2\r):\,\psi_1\c\s=\psi_1,\,\psi_2\c\s=\pm\psi_2\r\}$$
The same argument also works for $\rho_m$ and $\tau_m$ as well. This proves that $\mcal T(\mu,\cd)$ maps $\mcal X_{n,m}$ into $\mcal Y_{n,m}$.\\
We are left to show that $\ker\l(\pa_\phi\mcal T(\mu_n,0)\r)\cap\mcal X_{n,m}=\mrm{span}\l\{\l(\barr{c}0\\P_n^m(z)\cos(m\t)\earr\r)\r\}$.\\
By the property that $P_n^m(-z)=(-1)^{n+m}P_n^m(z)$ we get that the former element belongs to $\mcal X_{n,m}$. Hence the claim follows if we rule out all the other generators of the kernel.\vskip0.3cm
{\bf Step 1: if $\l(\barr{c}w_1(\t,z)\\w_2(\t,z)\earr\r)\in span\l\{\l(\barr{c}z\\0\earr\r),\l(\barr{c}\sqrt{1-z^2}\cos\t\\0\earr\r),\l(\barr{c}\sqrt{1-z^2}\sin\t\\0\earr\r)\r\}$ then $\l(\barr{c}w_1(\t,z)\\w_2(\t,z)\earr\r)\not\in\mcal X_{n,m}$.
}\\
First of all, there are no constant multiples of $\l(\barr{c}z\\0\earr\r)$, since the invariance with respect to $\s$ would imply
$$Cz=C\s(z)=-Cz,$$
hence $C=0$.\\
Now, assume that $\mcal X^{m,n}$ contains an element of the kind $\l(\barr{c}w_1\\0\earr\r)$, with $w_1(\t,z)=\sqrt{1-z^2}(A\cos\t+B\sin\t)$ for some real $A,B$. Then, we would have $w_1=w_1\c\rho_m$ and  the latter can be expressed by
$$w_1(\rho_m(\t,z))=\l(\cos\l(\fr{\pi}m\r)A+\sin\l(\fr{\pi}m\r)B\r)\sqrt{1-z^2}\cos\t+\l(-\sin\l(\fr{\pi}m\r)A+\cos\l(\fr{\pi}m\r)B\r)\sqrt{1-z^2}\sin\t.$$
Therefore, the pair $(A,B)$ must solve the linear system
\bequ
\begin{split}
\l\{\bll A=\cos\l(\fr{\pi}m\r)A+\sin\l(\fr{\pi}m\r)B\\
\\
B=-\sin\l(\fr{\pi}m\r)A+\cos\l(\fr{\pi}m\r)B\earr\r.
\end{split}
\eequ
which has no non-trivial solutions for any $m\geq 1$.\vskip0.2cm
{\bf Step 2: if $\l(\barr{c}w_1(\t,z)\\w_2(\t,z)\earr\r)\in span\l\{\l(\barr{c}0\\P_n^0(z)\earr\r),
\l(\barr{c}0\\P_n^l(z)\cos(m\t)\earr\r),\l(\barr{c}0\\P_n^l(z)\sin(m\t)\earr\r)\r\}_{l=1,\ds,n}$ then $\l(\barr{c}w_1(\t,z)\\w_2(\t,z)\earr\r)=\l(\barr{c}0\\P_n^m(z)\cos(m\t)\earr\r)$.
}\vskip0.2cm
First we exclude functions of the kind $\l(\barr{c}0\\w_2\earr\r)$ with $w_2(\t,z)=P_n^l(z)(A\cos(l\t)+B\sin(l\t))$ for $l=1,\ds,n,\,l\ne m$. Since we must have $w_2=-w_2\c\rho_m$  it holds
$$w_2(\rho_m(\t,z))=P_n^l(z)\l(\l(\cos\l(\fr{l}m\pi\r)A+\sin\l(\fr{l}m\pi\r)B\r)\cos(l\t)+\l(-\sin\l(\fr{l}m\pi\r)A+\cos\l(\fr{l}m\pi\r)B\r)\sin(l\t)\r).$$
Therefore, $A$ and $B$ now solve
$$\l\{\bll A=-\cos\l(\fr{l}m\pi\r)A-\sin\l(\fr{l}m\pi\r)B\\
B=\sin\l(\fr{l}m\pi\r)A-\cos\l(\fr{l}m\pi\r)B\earr\r.,$$
and we have nontrivial solutions if and if only $\fr{l}m$ is an odd integer. Since $l\ne m$ and $\fr{l}m\le\fr{n}m<3$ we deduce that $A=B=0$.\\
This argument also works, setting $l=B=0$, to rule out  $\l(\barr{c}0\\P_n^0\l(z\r)\earr\r)$.\\
Finally, we are left with $w_2(x)=P_n^m(z)(A\cos(m\t)+B\sin(m\t))$. Here using $\tau_m$ in place of $\rho_m$ we  have
$$P_n^m(z)(A\cos(m\t)+B\sin(m\t))=w_2(x)=-w_2(\tau_m(x))=P_n^m(z)(A\cos(m\t)-B\sin(m\t)),$$
which implies $B=0$. So the claim follows.
\epf\

\brem\label{rotaz}${}$\\
For any $\ph\in\S^1$ one can construct similar subspaces $\mcal X_{n,m,\ph}\sub W^{2,2}\l(\S^2\r)\x W^{2,2}\l(\S^2\r)$ and $\mcal Y_{n,m,\ph}\sub L^2\l(\S^2\r)\x L^2\l(\S^2\r)$ such that $\mcal T\l(\mcal X_{n,m,\ph}\r)\sub\mcal Y_{n,m,\ph}$ and the kernel of its derivative in $(\mu_n,0)$ is generated by $\l(\barr{c}0\\P_n^m(z)\cos(m\t+\ph)\earr\r)$. This is consistent with the invariance with respect to rotation of the operator $\mcal T$.\\
To prove this fact, just argue as in Lemma \ref{simm} replacing the reflection $\tau_m$ with $\tau_{m,\ph}:\t\to-\t+\fr{\pi-2\ph}m$.\\
In particular, for $\ph=-\fr{\pi}2$, we get $P_n^m(z)\sin(m\t)$. Notice that, in general, 
$$P_n^m(z)\cos(m\t+\ph)=\cos(\ph)P_n^m(z)\cos(m\t)-\sin(\ph)P_n^m(z)\sin(m\t)$$
is a combination of the two generators of the kernel.\\
Of course, since system \eqref{sistema} is rotationally invariant, in general this argument does not give rise to new solutions different from the previous ones.\\

More generally, we may consider a closed group $G\le O(3)$, a group omomorphism $\si:G\to\{-1,1\}$ and subspaces
\beqy
\mcal X_G&:=&\l\{(\phi_1,\phi_2)\in W^{2,2}\l(\S^2\r)\x W^{2,2}\l(\S^2\r);\,\phi_1\c g=\phi_1,\,\phi_2\c g=\si(g)\phi_2,\,\fa\,g\in G\r\};\\
\mcal Y_G&:=&\l\{(\psi_1,\psi_2)\in L^2\l(\S^2\r)\x L^2\l(\S^2\r);\,\psi_1\c g=\psi_1,\,\psi_2\c g=\si(g)\psi_2,\,\fa\,g\in G\r\}.
\eeqy
By the symmetry properties of $\mcal T$ we have $\mcal T(\mcal X_G)\sub\mcal Y_G$. Therefore, if we can find $G,\si$ such that $\ker\l(\pa_\phi\mcal T(\mu_n,0)\r)\cap\mcal X_G$ is $1$-dimensional and satisfies the transversality condition, then the same argument can be used to find solutions in $\mcal X_G$.\\
The spaces $\mcal X_{n,m,\ph},\,\mcal Y_{n,m,\ph}$ are obtained taking
$$G=\l\langle\s,\,\rho_{\fr{\pi}m},\,\tau_\fr{\pi-2\ph}m\r\rangle\q\q\q\q\q\q\si(\s)=(-1)^{n,m},\,\si\l(\rho_{\fr{\pi}m}\r)=\si\l(\tau_\fr{\pi-2\ph}m\r)=-1.$$
By choosing instead $G=O(2)$, namely the group of isometries fixing the $z$ axes, and $\si\eq1$ we recover the space of functions depending only on the $z$ variable. Such functions correspond to radially symmetric functions on the plane, therefore in this case we recover the radial branch found in \cite{ggw}.
\erem\

\blem\label{ran}${}$\\
The range of the operator $\pa_\phi\mcal T_{n,m}(\mu_n,0,0)$ has co-dimension $1$ and it is given by
$$\mrm R(\pa_\phi\mcal T_{n,m}(\mu_n,0,0))=\l\{(\psi_1,\psi_2)\in L^2\l(\S^2\r)\x L^2\l(\S^2\r):\,\int_{\S^2}\psi_2(\t,z)P_n^m(z)\cos(m\t)\mrm d\t\mrm dz=0\r\}$$
\elem\

\bpf${}$\\
By the definition of $\mcal T$ and the well-known spectral properties of $-\D_{\S^2}$, the range of $\pa_\phi\mcal T_{n,m}(\mu_n,0,0)$ coincides with the orthogonal of its kernel.\\
In Lemma \ref{simm} we showed that the kernel is $1$-dimensional and spanned by $\l(\barr{c}0\\P_n^m(z)\cos(m\t)\earr\r)$. From this one can easily derive the characterization of the range and this proves the claim.
\epf\

\section{Proof of the main theorem}\label{s4}

\bpf[Proof of Theorem \ref{teo}]${}$\\
We apply Theorem \ref{cr} to the operator $\mcal T_{n,m}:(-2,2)\x\mcal X_{n,m}\to\mcal Y_{n,m}$. We have to verify that all the hypotheses are satisfied.\\
It is immediate to see that $\mcal T_{n,m}(\mu,0,0)=(0,0)$ for all $\mu$ and that the first and second derivatives exist and are continuous. Lemma \ref{simm} ensures that $\ker\l(\pa_\phi\mcal T_{n,m}(\mu_n,0,0)\r)$ is $1$-dimensional and, by Lemma \ref{ran}, the range has also co-dimension $1$. Finally, the transversality condition is ensured by Lemma \ref{trasv}.\\
Therefore, from Theorem \ref{cr} there exists a branch of non-trivial solutions $(\mu_n(\e),\phi_{1,\e},\phi_{2,\e})$.\\
Since the bifurcation evolves in the direction $\l(\barr{c}0\\P_n^m(z)\cos(m\t)\earr\r)$, the functions $\phi_{i,\e}$ will satisfy
$$\l(\barr{c}\phi_{1,\e}(\t,z)\\\phi_{2,\e}(\t,z)\earr\r)=2\e\l(\barr{c}0\\P_n^m(z)\cos(m\t)\earr\r)+2\e\l(\barr{c}z_{1,\e}(\t,z)\\z_{2,\e}(\t,z)\earr\r)$$
with $z_{i,\e}$ satisfying $z_{i,0}=0$. Moreover, since $(\phi_{1,\e},\phi_{2,\e})\in\mcal X_{n,m}$, then also $(z_{1,\e},z_{2,\e})$ belongs to $\mcal X_{n,m}$ so it satisfies the symmetry properties:
$$\barr{lllll}
z_{1,\e}(\t,-z)=z_{1,\e}(\t,z)&\q\q&z_{1,\e}\l(-\t+\fr{\pi}m,z\r)=z_{1,\e}(\t,z)&\q\q&z_{1,\e}\l(\t+\fr{\pi}m,z\r)=z_{1,\e}(\t,z)\\
z_{2,\e}(\t,-z)=(-1)^{n+m}z_{2,\e}(\t,z)&\q\q&z_{2,\e}\l(-\t+\fr{\pi}m,z\r)=-z_{2,\e}(\t,z)&\q\q&z_{2,\e}\l(\t+\fr{\pi}m,z\r)=-z_{2,\e}(\t,z)\\
\earr.$$
By Proposition \ref{sfera} we have that the functions $(u_{1,\e},u_{2,\e})$ defined as in \eqref{u} solve \eqref{sistema} and they satisfy all the requested properties. In particular, the symmetries \eqref{z} satisfied by $(Z_{1,\e},Z_{2,\e})=\l(z_{1,\e}\c\Pi^{-1},z_{2,\e}\c\Pi^{-1}\r)$ follow from Remark \ref{simmetrie}.\\
Moreover since $z_{1,\e},z_{2,\e}\in W^{2,2}(\S^2)$ then they are bounded in $L^\infty$ and the same happens for $(Z_{1,\e},Z_{2,\e})$. Hence we get that $(u_{1,\e},u_{2,\e})$  verify
$$\l\{\bl u_{1,\e}(x)=U_{\mu_n(\e)}(x)+\e P_n^m\l(\fr{8-|x|^2}{8+|x|^2}\r)\cos(m\t)+\e(Z_{1,\e}(x)+Z_{2,\e}(x))\\
u_{2,\e}(x)=U_{\mu_n(\e)}(x)-\e P_n^m\l(\fr{8-|x|^2}{8+|x|^2}\r)\cos(m\t)+\e(Z_{1,\e}(x)-Z_{2,\e}(x)),\earr\r.$$
In order to get \eqref{i20} observe that using the Taylor expansion we get
\bequ\label{1}
U_{\mu_n(\e)}(x)=U_{\mu_n}(x)+\log\frac1{2+\mu_n(\e)}-\log\frac1{2+\mu_n}=
U_{\mu_n}(x)-\frac{\mu'_n(0)}{2+\mu_n}\e+C(\e).
\eequ
where $C(\e)$ is a term satisfying  $C(\e)=O(\e^2)$. In Lemma \ref{mu0} in Appendix we prove that $\mu_n'(0)=0$ and so the claim follows replacing $Z_{1,\e}$ with $Z_{1,\e}+\frac{C(\e)}\e$.\\
Crandall-Rabinowitz Theorem provides  the existence of {\em local}  branches. In order to prove that the bifurcation is global we use a degree argument (see \cite{kie}). Since we need to define an operator from a Banach space into itself let us consider $\mcal L:L^\infty(\S^2)\rightarrow L^\infty(\S^2)$ as
\bequ\label{l}
\mcal L:(\mu,\phi_1,\phi_2)=\l(\barr{c}
\phi_1-(\D_{\S^2}+I)^{-1}\l[2\l(e^\fr{\phi_1+\phi_2}2+e^\fr{\phi_1-\phi_2}2-2\r)+\phi_1\r]\\ \\
\phi_2-(\D_{\S^2}+I)^{-1}\l[2\fr{2-\mu}{2+\mu}\l(e^\fr{\phi_1+\phi_2}2-e^\fr{\phi_1-\phi_2}2\r)+\phi_2\r]\earr\r).\eequ
It is immediate to check that Lemma \ref{simm} holds for the operator $\mcal L$ in the space $\mcal X_{n,m}$. This implies that as the parameter $\mu$ crosses $\mu_n$ the Morse index of the solution increases by {\em one}. Hence there is a change in the Leray-Schauder degree $\deg(\mcal L,B_\delta,0)$ and then classical results in bifurcation theory ensure
 the existence of a global continuum of solutions to \eqref{b1} that satisfies the Rabinowitz alternative.\\
Now let us prove \eqref{i4}. Integrating \eqref{b1} we get that
$$\int_{\S^2}e^{\frac{\phi_{1,\e}+\phi_{2,\e}}2}=\int_{\S^2}e^{\frac{\phi_{1,\e}-\phi_{2,\e}}2}=4\pi.
$$
Then we get
$$\int_{\R^2}e^{u_{1,\e}}=\int_{\R^2}e^{U_{\mu_n(\e)}}e^\fr{\Phi_{1,\e}+\Phi_{2,\e}}2=\fr{2}{2+\mu_n(\e)}\int_{\S^2}e^{\frac{\phi_{1,\e}+\phi_{2,\e}}2}=\fr{8\pi}{2+\mu_n(\e)}$$
and the same holds for $\int_{\R^2}e^{u_{2,\e}}$. This shows  \eqref{i4}.
\epf\

\bpf[Proof of Corollary \ref{cor2}]${}$\\
By Theorem \ref{teo} we get the existence of $n-\left[\frac n3\right]$ nonradial nonequivalent solutions and in \cite{ggw} we got the existence of a radial solution bifurcating by $(\mu_n,U_{\mu_n},U_{\mu_n})$. So the claim follows.
\epf\

\brem\label{e1}
Let us consider the case $n=2$, namely the very well-known $SU(3)$ Toda system.\\
Here, the assumption $m>\fr{n}3$ is redundant and Theorem \ref{teo} can be applied for $m=1,2$ to get two branches of non-trivial solutions.\\
Through Remark \ref{rotaz} we get two more independent branches of non-radial solutions; moreover, in \cite{ggw}, Theorem $1.2$ the authors find a non-trivial branch of radial solutions.\\
Therefore, we get five different non-trivial branches of solution. Putting them together with the $3$-parameters family of \qm{trivial} solutions $\{U_{\mu_2,\d,y}\}_{\d\in\R,y\in\R^2}$ we recover  the $8$-parameters family of explicit solutions found in \cite{jw,lwy}.
\erem\

We end this section giving the proof of Corollary \ref{cor3}.
\bpf[Proof of Corollary \ref{cor3}]${}$\\ 
 We consider the operator $\mcal L:L^\infty(\S^2)\rightarrow L^\infty(\S^2)$
 as defined in \eqref{l} restricted to space of functions which depend only on the variable $z$. Using \eqref{lp} a zero of $\mcal L$ then solves
$$\l\{\bl 
-(1-z^2) \partial_{zz}\phi_1(z)+2z\partial_z \phi_1   =2\left( e^{\fr {\phi_1+\phi_2}2}+e^{\fr {\phi_1-\phi_2}2}-2\right)\\
-(1-z^2) \partial_{zz}\phi_2(z)+2z\partial_z \phi_2   =2 \fr   {2-\mu}{2+\mu}\left( e^{\fr {\phi_1+\phi_2}2}-e^{\fr {\phi_1-\phi_2}2}\right),
\earr\r.$$ 
 This is a Sturm Liouville problem, and this implies that the solutions can only have simple zeroes. \\
 Near the bifurcation point $\mu_n$ we know, from the local bifurcation result, that   
 $$\phi_{2,\e}(z)=2\e P_n^0(z)+2\e z_{2,\e}(z)$$
 and it is not hard to prove that $ \phi_{2,\e}(z)$ has the same number of zeroes of $P_n(z)$, since $P_n(z)$ has only simple zeroes. Then, by standard Sturm Liouville theory, the number of zeroes of $ \phi_{2,\e}(z)$ cannot change along the branch. Since $P_n(z)$ has $n$ zeroes in $(-1,1)$ this implies that the radial branch bifurcating from $\mu_n$ cannot join the radial branch bifurcating from $\mu_l$ for $n\neq l$ and so they do not intersect. 
\epf\

\appendix

\section{Appendix. The derivatives of the parameter $\mu_n(\e)$}\label{s6}
In this appendix we study the behavior of the parameter $\mu_n(\e)$ in a neighborhood of $\mu_{n}$. Let us start recalling a known result in bifurcation theory.
\bthm{\cite{kie}, Chapter $I.6$}${}$\\
Assume the hypotheses of Theorem \ref{cr} are satisfied and, moreover, $F\in C^2((-2,2)\x X,Y)$.\\
Let $y_0\in Y$ be any generator of $\mrm R(\pa_xF(0,0))^\perp$. Then, the derivative $\mu'(0)$ of $\mu(\e)$ in $\e=0$ is given by
\bequ\label{mu1}
\mu'(0)=-\fr{1}2\fr{\l\langle\pa^2_{x,x}F(0,0)[w_0,w_0],y_0\r\rangle}{\|w_0\|\l\langle\pa^2_{t,x}F(0,0)w_0,y_0\r\rangle}.
\eequ
Assume, in addition, that $\mu'(0)=0$ and $F$ is of class $C^3$. Define as $\mrm Q:y\to\fr{\langle y,y_0\rangle}{\|y_0\|^2}y_0$ the projection from $Y$ to $\mrm R(\pa_xF(0,0))^\perp$ and write as $\l(\pa_xF(0,0)\r)^{-1}:\mrm R(\pa_xF(0,0))\to\ker(\pa_xF(0,0))^\perp$ the inverse of $\pa_xF(0,0)$ restricted to the complementary of its kernel.\\
Then, $\mu''(0)$ is given by
\bequ\label{mu2}
\mu''(0)=-\fr{1}3\fr{\l\langle\pa^3_{x,x,x}F(0,0)[w_0,w_0,w_0]-3\pa^2_{x,x}F(0,0)\l[w_0,\l(\pa_xF(0,0)\r)^{-1}(\mrm I-\mrm Q)\pa^2_{x,x}F(0,0)[w_0,w_0]\r],y_0\r\rangle}{\|w_0\|^2\l\langle\pa^2_{t,x}F(0,0)w_0,y_0\r\rangle}.
\eequ
\ethm\

\blem\label{mu0}${}$\\
For any $n,m\in\N$ with $1\le\fr{n}m<3$, the parameter $\mu_n(\e)$ defined in Theorem \ref{teo} satisfies $\mu'_n(0)=0$.
\elem
\bpf${}$\\
We suffice to apply the formula \eqref{mu1} of Theorem \ref{teo}.\\
We recall that, in our case, $w_0=\l(\barr{c}0\\P_n^m(z)\cos(m\t)\earr\r)$ and the second derivative of $\mcal T_{n,m}$ in $0$ is
\bequ\label{d2}
\pa^2_{\phi,\phi}\mcal T_{n,m}(\mu_n,0,0):\l(\l(\barr{c}v_1\\v_2\earr\r),\l(\barr{c}w_1\\w_2\earr\r)\r)\to\l(\barr{c}v_1w_1+v_2w_2\\\fr{n(n+1)}2(v_1w_2+v_2w_1)\earr\r).
\eequ
Moreover, thanks to Lemma \ref{ran}, the orthogonal of the range is generated by $y_0=w_0$.\\
Therefore, \eqref{mu1} gives:
\beqy
\mu'(0)&=&-\fr{1}2\fr{\l\langle\pa^2_{x,x}F(0,0)[w_0,w_0],w_0\r\rangle}{\|w_0\|\l\langle\pa^2_{t,x}F(0,0)w_0,w_0\r\rangle}\\
&&-\fr{\l\langle\pa^2_{\phi,\phi}\mcal T_{n,m}(\mu_n,0,0)\l[\l(\barr{c}0\\P_n^m(z)\cos(m\t)\earr\r),\l(\barr{c}0\\P_n^m(z)\cos(m\t)\earr\r)\r],\l(\barr{c}0\\P_n^m(z)\cos(m\t)\earr\r)\r\rangle}{2\|w_0\|\l\langle\pa^2_{\mu,\phi}\mcal T_{n,m}(\mu_n,0,0)w_0,w_0\r\rangle}\\
&=&-\fr{\l\langle\l(\barr{c}\l(P_n^m(z)\cos(m\t)\r)^2\\0\earr\r),\l(\barr{c}0\\P_n^m(z)\cos(m\t)\earr\r)\r\rangle}{2\|w_0\|\l\langle\pa^2_{\mu,\phi}\mcal T_{n,m}(\mu_n,0,0)w_0,w_0\r\rangle}\\
&=&0
\eeqy
\epf
The following functions play a similar role with respect to the Legendre polynomials \eqref{i3}:
\bequ\label{pnmtilde}
\wt P_n^m(z)=\l(\fr{1-z}{1+z}\r)^\fr{m}2\sum_{k=0}^n(-1)^k\bin{n}k\fr{\bin{k+n}k}{\bin{k+m}k}\l(\fr{1-z}2\r)^k
\eequ
for $n,m\in \N$.

\blem\label{soluz}${}$\\
For any $\l(\barr{c}\psi_1\\\psi_2\earr\r)\in\mrm R(\pa_\phi\mcal T(\mu_n,0,0))^\perp$, the only solution $\l(\barr{c}\phi_1\\\phi_2\earr\r)\in\ker\l(\pa_\phi\mcal T(\mu_n,0,0)\r)^\perp$ of
$$\pa_\phi\mcal T(\mu_n,0,0)\l(\barr{c}\phi_1\\\phi_2\earr\r)=\l(\barr{c}\D_{\S^2}\phi_1+2\phi_1\\
\D_{\S^2}\phi_2+n(n+1)\phi_2\earr\r)=\l(\barr{c}\psi_1\\\psi_2\earr\r)$$
is given by
$$\phi_i(\t,z):=\phi_{i,0}(z)+\sum_{l=1}^{+\infty}\l(\phi_{i,l}^1(z)\cos(m\t)+\phi_{i,l}^2(z)\sin(m\t)\r),$$
where
\beqa
\nonumber\phi_{1,0}(z)&=&P_1(z)\l(C_{1,0}-\int_{-1}^z\fr{1}{\l(1-y^2\r)(P_1(y))^2}\int_y^1P_1(x)\psi_{1,0}(x)\mrm dx\mrm dy\r)\\
\nonumber\phi_{1,1}^j(z)&=&P_1^1(z)\l(C_{1,1}^j-\int_0^z\fr{1}{\l(1-y^2\r)(P_1^1(y))^2}\int_y^1P_1^1(x)\psi_{1,1}^j(x)\mrm dx\mrm dy\r)\\
\nonumber(l\ge2)\q\phi_{1,l}^j(z)&=&-\wt P_1^l(z)\int_{-1}^z\fr{1}{\l(1-y^2\r)\l(\wt P_1^l(y)\r)^2}\int_y^1\wt P_1^l(x)\psi_{1,l}^j(x)\mrm dx\mrm dy\\
\label{phi20}\phi_{2,0}(z)&=&P_n(z)\l(C_{2,0}-\int_0^z\fr{1}{\l(1-y^2\r)(P_n(y))^2}\int_y^1P_n(x)\psi_{n,0}(x)\mrm dx\mrm dy\r)\\
\nonumber(l\le n)\q\phi_{2,l}^j(z)&=&P_n^l(z)\l(C_{2,l}^j-\int_{-1}^z\fr{1}{\l(1-y^2\r)(P_n^l(y))^2}\int_y^1P_n^l(x)\psi_{2,l}^j(x)\mrm dx\mrm dy\r)\\
\nonumber(l\ge n+1)\q\phi_{2,l}^j(z)&=&-\wt P_n^l(z)\int_{-1}^z\fr{1}{\l(1-y^2\r)\l(\wt P_n^l(y)\r)^2}\int_y^1\wt P_n^l(x)\psi_{2,l}^j(x)\mrm dx\mrm dy\\
\eeqa
with $C_{i,0},C_{i,l}^j$ uniquely determined and $\psi_{i,0},\psi_{i,l}^j$ defined by the Fourier decomposition in $\t$:
$$\psi_i(\t,z)=\psi_{i,0}(z)+\sum_{l=1}^{+\infty}\l(\psi_{i,l}^1(z)\cos(l\t)+\psi_{i,l}^2(z)\sin(l\t)\r)$$
\elem\

\bpf${}$\\
Since the linearized operator has two decoupled equations, we will suffice to solve $\D\phi+n(n+1)\phi=\psi$ for any $\psi$ being orthogonal to the eigenfunctions of the homogeneous problem, which we recall are given by \eqref{base}.\\
In the system of coordinates $(\t,z)$ the latter problem has the form
$$\l(1-z^2\r)\pa_{zz}\phi(\t,z)-2z\pa_{z}\phi(\t,z)+\fr{1}{1-z^2}\pa_{\t\t}\phi(\t,z)+n(n+1)\phi(\t,z)=\psi(\t,z).$$
If we now decompose $\phi,\psi$ as
\beqy
\phi(\t,z)&=&\phi_0(z)+\sum_{l=1}^{+\infty}\l(\phi_l^1(z)\cos(l\t)+\phi_l^2(z)\sin(l\t)\r)\\
\psi(\t,z)&=&\psi_0(z)+\sum_{l=1}^{+\infty}\l(\psi_l^1(z)\cos(l\t)+\psi_l^2(z)\sin(l\t)\r),
\eeqy
then each $\phi_l^j$ solves the o.d.e.
$$\l(1-z^2\r){\phi_l^j}''(z)-2z{\phi_l^j}'(z)-\fr{l^2}{1-z^2}\phi_l^j(z)+n(n+1)\phi_l^j(z)=\psi_l^j(z).$$
If $l\le n$ and $\psi_l^j(z)\equiv0$, the only bounded solution of this equation, up to a multiplicative constant, is $P_n^l(z)$; on the other hand, for $l\ge n+1$, we just get a solution $\wt P_n^l(z)$ defined by \eqref{pnmtilde} which is bounded in $\S^2\sm\{0,0,-1\}$.\\
Now, using the variation of the constants we get
\bequ\label{varcost}
\phi_l^j(z)=P_n^l(z)\int_a^z\fr{1}{\l(1-y^2\r)(P_n^l(y))^2}\int_b^yP_n^l(x)\psi_l^j(x)\mrm dx\mrm dy
\eequ
for $l\le n$ and a similar expression, with $\wt P_n^l$ in place of $P_n^l$. Here $a$ and $b$ are to be determined in such a way that the function is globally defined on the whole sphere, namely up to $y=\pm1$\\
It is not hard to see that to get boundedness in $z=\pm1$ we must have  $a=-1$ and $b=1$.\\
Arguing similarly we get
\bequ\label{c}
C_0:=\fr{\int_{-1}^1(P_n(z))^2\int_{-1}^z\fr{1}{\l(1-y^2\r)(P_n(y))^2}\int_y^1P_n(x)\psi_0(x)\mrm dx\mrm dy\mrm dz}{\int_{-1}^1(P_n^l(z))^2\mrm dz}\q\q\q C_l^j:=\fr{\int_{-1}^1(P_n^l(z))^2\int_0^z\fr{1}{\l(1-y^2\r)(P_n^l(y))^2}\int_y^1P_n^l(x)\psi_l^j(x)\mrm dx\mrm dy\mrm dz}{\int_{-1}^1(P_n^l(z))^2\mrm dz}.
\eequ
By putting together the formulas \eqref{varcost} and \eqref{c} we conclude the proof.
\epf\

\brem${}$\\
We will actually need to invert not $\pa_\phi\mcal T(\mu_n,0,0)$ but just $\pa_\phi\mcal T_{n,m}(\mu_n,0,0)$, namely the restriction of $\pa_\phi\mcal T(\mu_n,0,0)$ on the spaces $\mcal X_{n,m},\mcal Y_{n,m}$ defined by \eqref{xnm},\eqref{ynm}.\\
Choosing $(\psi_1,\psi_2)\in\mcal Y_{n,m}$ some of the coefficients $\psi_{i,l}^j$ will vanish (hence also $\phi_{i,l}^j$ will). In particular, arguing as in the proof of Lemma \ref{simm} one can see that $\psi_{i,l}^2=0$ for any $i,l$, $\psi_{1,l}^1=0$ unless $l$ is an even multiple of $m$ and $\psi_{1,l}^2=0$ unless $l$ is an odd multiple of $m$ (and all the non-vanishing $\psi_{i,0},\psi_{i,l}^j$ are even or odd in $z$).\\
Anyway, the argument to invert the linearized operator of $\mcal T_{n,m}$ still works for the more general $\mcal T$ and for this reason we gave the proof for the more general case.
\erem\
\begin{remark}${}$\\
The formulas for $\phi_{i,0},\phi_{i,l}^j$ in Lemma \ref{soluz} are not really well-defined for all $z$, because the denominators of the fractions may vanish. In fact, the integral \eqref{phi20} defining $\phi_{2,0}$ makes no sense for $z$ greater than the first zero $z_1$ of $P_n(z)$ and similar issues occur in the other formulas.\\
To fix this, we have to replace the integral with some anti-derivatives of the same function and choosing properly the free constant.\\
Since $P_n$ has only simple zeros $z_1<\ds<z_n$, then the limit of \eqref{phi20} as $z$ goes to $z_1$ makes sense. Then, we can define
$$\phi_{2,0}=P_n(z)\l(C_{2,0}-\int_{a_n}^z\fr{\mrm dy}{\l(1-y^2\r)(P_n(y))^2}\int_y^1P_n(x)\psi_{2,0}(x)\mrm dx\r)\q\q\q\tx{if }z_1<z<z_2$$
with $a_n\in(z_1,z_2)$. This definition is also well-posed as $z$ goes to $z_1$ and, for a proper choice of $a_n$, the two definitions can be \qm{glued} smoothly around $z=z_1$. We then repeat the procedure for each interval $(z_2,z_3),\ds,(z_n,1)$ and we find a smooth solution for all $z$.\\
The same argument also fixes other formulas introduced in Lemma \ref{soluz}, since the zeros of both $P_n^l$ and $\wt P_n^l$ are also simple.
\end{remark}
\blem${}$\\
For any $n,m\in\N$ with $1\le\fr{n}m<3$ the parameter $\mu_n(\e)$ defined in Theorem \ref{teo} satisfies
\beqy\label{h}
\mu''(0)&=&C_{m,n}\l(\int_0^1(P_n^m(z))^4\mrm dz+2\int_{-1}^1z(P_n^m(z))^2\int_{-1}^z\fr{1}{y^2\l(1-y^2\r)}\int_y^1x(P_n^m(x))^2\mrm dx\mrm dy\mrm dz\r.\\
&+&\l.\int_{-1}^1(z+2m)\l(\fr{1-z}{1+z}\r)^m(P_n^m(z))^2\int_{-1}^z\fr{1}{(y+2m)^2\l(1-y^2\r)\l(\fr{1-y}{1+y}\r)^{2m}}\int_y^1(x+2m)\l(\fr{1-x}{1+x}\r)^m(P_n^m(x))^2\mrm dx\mrm dy\mrm dz\r).
\eeqy
with $C_{n,m}:=\fr{n(n+1)}{4\l(n^2+n+1\r)\pi}\l(\fr{(2n+1)(n-m)!}{\l(n^2+n+2\r)(n+m)!}\r)^2>0$.
\elem\

\bpf${}$\\
Since, by Lemma \ref{mu0}, we have $\mu'(0)=0$, then we are in position to apply the formula \eqref{mu2} to evaluate $\mu''(0)$.\\
To this purpose, we recall the second and third derivative of $\mcal T$, given respectively by \eqref{d2} and
$$\pa^3_{\phi,\phi,\phi}\mcal T_{m,n}(\mu_n,0,0):\l(\l(\barr{c}u_1\\u_2\earr\r),\l(\barr{c}v_1\\v_2\earr\r),\l(\barr{c}w_1\\w_2\earr\r)\r)\to\l(\barr{c}\fr{u_1v_1w_1+u_1v_2w_2+u_2v_1w_2+u_2v_2w_1}2\\n(n+1)\fr{u_1v_1w_2+u_1v_2w_1+u_2v_1w_1+u_2v_2w_2}4\earr\r).$$
We also need the explicit expression of the inverse operator $\l(\pa_\phi\mcal T_{m,n}(\mu_n,0,0)\r)^{-1}:\mrm R\l(\pa_\phi\mcal T(\mu_n,0,0)\r)\to\ker\l(\pa_\phi\mcal T_{m,n}(\mu_n,0,0)\r)^\perp$, which is given by Lemma \eqref{soluz}.\\
Concerning the denominator in \eqref{mu2}, as in the proof of Lemma \ref{trasv} we have
\beqy
\|w_0\|^2\l\langle\pa^2_{\mu,\phi}\mcal T_{m,n}(\mu_n,0,0)w_0,y_0\r\rangle&=&\l(\int_{\S^2}\l(|\n w_0|^2+|w_0|^2\r)\r)\l\langle\l(\barr{c}0\\-\fr{8}{(2+\mu_n)^2}P_n^m(z)\cos(m\t)\earr\r),\l(\barr{c}0\\P_n^m(z)\cos(m\t)\earr\r)\r\rangle\\
&=&\l(\int_{\S^2}w_0\cd(-\D_{\S^2} w_0+w_0)\r)\l(-\fr{8}{(2+\mu_n)^2}\int_{\S^2}\l(P_n^m(z)\cos(m\t)\r)^2\mrm d\t\mrm dz\r)\\
&=&-\fr{8}{(2+\mu_n)^2}(n(n+1)+1)\l(\int_{\S^2}|w_0|^2\r)\l(\int_0^{2\pi}(\cos(m\t))^2\mrm d\t\int_{-1}^1(P_n^m(z))^2\mrm dz\r)\\
&=&-\fr{\l(n^2+n+2\r)^2\l(n^2+n+1\r)}8\l(\int_0^{2\pi}(\cos(m\t))^2\mrm d\t\int_{-1}^1(P_n^m(z))^2\mrm dz\r)^2\\
&=&-\fr{\l(n^2+n+2\r)^2\l(n^2+n+1\r)}8\l(\pi\fr{2(n+m)!}{(2n+1)(n-m)!}\r)^2\\
\eeqy
The first term in the numerator of \eqref{mu2} is not hard to compute:
\beqy
&&\l\langle\pa^3_{\phi,\phi,\phi}\mcal T(\mu_n,0,0)\l[\l(\barr{c}0\\P_n^m(z)\cos(m\t)\earr\r)\l(\barr{c}0\\P_n^m(z)\cos(m\t)\earr\r)\l(\barr{c}0\\P_n^m(z)\cos(m\t)\earr\r)\r],\l(\barr{c}0\\P_n^m(z)\cos(m\t)\earr\r)\r\rangle\\
&=&\l\langle\l(\barr{c}0\\\fr{n(n+1)}4((P_n^m(z)\cos(m\t))^3\earr\r),\l(\barr{c}0\\P_n^m(z)\cos(m\t)\earr\r)\r\rangle\\
&=&\fr{n(n+1)}4\int_0^{2\pi}\cos(m\t)^4\mrm d\t\int_{-1}^1(P_n^m(z))^4\mrm dz\\
&=&\fr{3n(n+1)}{8}\pi\int_0^1(P_n^m(z))^4\mrm dz
\eeqy

To evaluate the other term in the numerator, we first notice that
\beqy
\mrm Q\pa^2_{\phi,\phi}\mcal T_{m,n}(\mu_n,0)[w_0,w_0]&=&\fr{\l\langle\pa^2_{\phi,\phi}\mcal T(\mu_n,0,0)\l[\l(\barr{c}0\\P_n^m(z)\cos(m\t)\earr\r),\l(\barr{c}0\\P_n^m(z)\cos(m\t)\earr\r)\r],y_0\r\rangle}{\|y_0\|^2}y_0\\
&=&\fr{1}{\|y_0\|^2}\l\langle\l(\barr{c}(P_n^m(z)\cos(m\t))^2\\0\earr\r)\l(\barr{c}0\\P_n^m(z)\cos(m\t)\earr\r)\r\rangle y_0\\
&=&0;
\eeqy
therefore, we get
\beqy
&&\l(\pa_\phi\mcal T(\mu_n,0,0)\r)^{-1}(\mrm I-\mrm Q)\pa^2_{\phi,\phi}\mcal T(\mu_n,0,0)\l[w_0,w_0\r]\\
&=&\l(\pa_\phi\mcal T(\mu_n,0,0)\r)^{-1}\pa^2_{\phi,\phi}\mcal T(\mu_n,0,0)\l[\l(\barr{c}0\\P_n^m(z)(\cos m\t)\earr\r),\l(\barr{c}0\\P_n^m(z)\cos(m\t)\earr\r)\r]\\
&=&\l(\pa_\phi\mcal T(\mu_n,0,0)\r)^{-1}\l(\barr{c}(P_n^m(z)\cos(m\t))^2\\0\earr\r)\\
&=&\l(\pa_\phi\mcal T(\mu_n,0,0)\r)^{-1}\l(\barr{c}(P_n^m(z))^2\l(\fr{1}2+\fr{\cos(2m\t)}2\r)\\0\earr\r).
\eeqy
We can apply Lemma \ref{soluz} to the latter. The only non-vanishing coefficients are $\psi_{1,0}(z)=\psi_{1,2m}^1(z)=\fr{1}2(P_n^m(z))^2$. Therefore, we get $\l(\barr{c}\phi_1(z)\\0\earr\r)$ with
\beqy
\phi_1(z)&=&P_1(z)\l(C_{1,0}-\fr{1}2\int_{-1}^z\fr{1}{\l(1-y^2\r)(P_1(y))^2}\int_y^1P_1(x)(P_n^m(x))^2\mrm dx\mrm dy\r)\\
&-&\fr{1}2\wt P_1^{2m}(z)\int_{-1}^z\fr{1}{\l(1-y^2\r)\l(\wt P_1^{2m}(y)\r)^2}\int_y^1\wt P_1^{2m}(x)(P_n^m(x))^2\mrm dx\mrm dy\cos(2m\t)
\eeqy
Therefore,
\beqy
&&\l\langle\pa^2_{\phi,\phi}\mcal T_{m,n}(\mu_n,0,0)\l[w_0,\l(\pa_\phi\mcal T_{m,n}(\mu_n,0,0)\r)^{-1}(\mrm I-\mrm Q)\pa^2_{\phi,\phi}\mcal T_{m,n}(\mu_n,0,0)[w_0,w_0]\r],y_0\r\rangle\\
&=&\l\langle\pa^2_{\phi,\phi}\mcal T_{m,n}(\mu_n,0,0)\l[\l(\barr{c}0\\P_n^m(z)\cos(m\t)\earr\r),\l(\barr{c}\phi_1(z,\t)\\0\earr\r)\r],y_0\r\rangle\\
&=&\l\langle\l(\barr{c}0\\\fr{n(n+1)}2P_n^m(z)\cos(m\t)\phi_1(z,\t)\earr\r),\l(\barr{c}0\\P_n^m(z)\cos(m\t)\earr\r)\r\rangle\\
&=&\fr{n(n+1)}2\int_{\S^2}(P_n^m(z)\cos(m\t))^2\phi_1(z,\t)\mrm d\t\mrm dz\\
&=&\fr{n(n+1)}2\int_0^{2\pi}\mrm d\t\int_{-1}^1\mrm dz(P_n^m(z)\cos(m\t))^2\l(P_1(z)\l(C_{1,0}-\fr{1}2\int_{-1}^z\fr{1}{\l(1-y^2\r)(P_1(y))^2}\int_y^1P_1(x)(P_n^m(x))^2\mrm dx\mrm dy\r)\r.\\
&-&\l.\fr{1}2\wt P_1^{2m}(z)\int_{-1}^z\fr{1}{\l(1-y^2\r)\l(\wt P_1^{2m}(y)\r)^2}\int_y^1\wt P_1^{2m}(x)(P_n^m(x))^2\mrm dx\mrm dy\cos(2m\t)\r)\\\\
&=&\fr{n(n+1)}2C_{1,0}\int_0^{2\pi}(\cos(m\t))^2\int_{-1}^1P_1(z)(P_n^m(z))^2\mrm dz\\
&-&\fr{n(n+1)}4\int_0^{2\pi}(\cos(m\t))^2\mrm d\t\int_{-1}^1P_1(z)(P_n^m(z))^2\int_{-1}^z\fr{1}{\l(1-y^2\r)(P_1(y))^2}\int_y^1P_1(x)(P_n^m(x))^2\mrm dx\mrm dy\mrm dz\\
&-&\fr{n(n+1)}4\int_0^{2\pi}(\cos(m\t))^2\cos(2m\t)\mrm d\t\int_{-1}^1\wt P_1^{2m}(z)(P_n^m(z))^2\int_{-1}^z\fr{1}{\l(1-y^2\r)\l(\wt P_1^{2m}(y)\r)^2}\int_y^1\wt P_1^{2m}(x)(P_n^m(x))^2\mrm dx\mrm dy\mrm dz\\
&=&-\fr{n(n+1)}4\pi\int_{-1}^1z(P_n^m(z))^2\int_{-1}^z\fr{1}{\l(1-y^2\r)y^2}\int_y^1x(P_n^m(x))^2\mrm dx\mrm dy\mrm dz\\
&-&\fr{n(n+1)}8\pi\int_{-1}^1(z+2m)\l(\fr{1-z}{1+z}\r)^m(P_n^m(z))^2\int_{-1}^z\fr{1}{\l(1-y^2\r)(y+2m)^2\l(\fr{1-y}{1+y}\r)^{2m}}\int_y^1(x+2m)\l(\fr{1-x}{1+x}\r)^m(P_n^m(x))^2\mrm dx\mrm dy\mrm dz\\
\eeqy
We are now in position to apply the formula \eqref{mu2} and to prove the Lemma:
\beqy
\mu''(0)&=&\fr{1}{\|w_0\|^2\l\langle\pa^2_{\mu,\phi}\mcal T_{n,m}(\mu_n,0,0)w_0,y_0\r\rangle}\l(-\fr{1}3\l\langle\pa^3_{\phi,\phi,\phi}\mcal T_{n,m}(\mu_n,0,0)[w_0,w_0,w_0],y_0\r\rangle\r.\\
&+&\l.\l\langle\pa^2_{\phi,\phi}\mcal T_{n,m}(\mu_n,0,0)\l[w_0,\l(\pa_\phi\mcal T_{m,n}(\mu_n,0,0)\r)^{-1}\l(\mrm I-\mrm Q\r)\pa^2_{\phi,\phi}\mcal T_{m,n}(\mu_n,0,0)[w_0,w_0]\r],y_0\r\rangle\r)\\
&=&-\fr{8}{\l(n^2+n+2\r)^2\l(n^2+n+1\r)}\l(\fr{(2n+1)(n-m)!}{2\pi(n+m)!}\r)^2\l(-\fr{1}3\fr{3n(n+1)}{8}\pi\int_0^1(P_n^m(z))^4\mrm dz\r.\\
&-&\l.\fr{n(n+1)}4\pi\int_{-1}^1z(P_n^m(z))^2\int_{-1}^z\fr{1}{\l(1-y^2\r)y^2}\int_y^1x(P_n^m(x))^2\mrm dx\mrm dy\mrm dz\r.\\
&-&\l.\fr{n(n+1)}8\pi\int_{-1}^1(z+2m)\l(\fr{1-z}{1+z}\r)^m(P_n^m(z))^2\int_{-1}^z\fr{1}{\l(1-y^2\r)(y+2m)^2\l(\fr{1-y}{1+y}\r)^{2m}}\int_y^1(x+2m)\l(\fr{1-x}{1+x}\r)^m(P_n^m(x))^2\mrm dx\mrm dy\mrm dz\r)\\
&=&C_{n,m}\l(\int_0^1(P_n^m(z))^4\mrm dz+2\int_{-1}^1z(P_n^m(z))^2\int_{-1}^z\fr{1}{\l(1-y^2\r)y^2}\int_y^1x(P_n^m(x))^2\mrm dx\mrm dy\mrm dz\r.\\
&+&\l.\int_{-1}^1(z+2m)\l(\fr{1-z}{1+z}\r)^m(P_n^m(z))^2\int_{-1}^z\fr{1}{\l(1-y^2\r)(y+2m)^2\l(\fr{1-y}{1+y}\r)^{2m}}\int_y^1(x+2m)\l(\fr{1-x}{1+x}\r)^m(P_n^m(x))^2\mrm dx\mrm dy\mrm dz\r)\\
\eeqy
\epf
\begin{remark}
We are not able to compute the integrals in \eqref{h} for any $m,n$. Of course this is possible in some particular cases like in the table below
$$
\barr{|l|l|l|l|}
\hline
n&\mu''(0)>0&\mu''(0)<0&\mu''(0)>0\\
\hline
3&m=0,1&m=2&m=3\\
\hline
4&m=0,1&m=2,3&m=4\\
\hline
5&m=0,1&m=2,3,4&m=5\\
\hline
6&m=0,1&m=2,3,4&m=5,6\\
\hline
7&m=0,1,2&m=3,4,5&m=6,7\\
\hline
8&m=0,1,2&m=3,4,5,6&m=7,8\\
\hline
9&m=0,1,2,3&m=4,5,6&m=7,8,9\\
\hline
10&m=0,1,2,3&m=4,5,6,7&m=8,9,10\\
\hline
11&m=0,1,2,3&m=4,5,6,7,8&m=9,10,11\\
\hline
12&m=0,1,2,3,4&m=5,6,7,8&m=9,10,11,12\\
\hline
13&m=0,1,2,3,4&m=5,6,7,8,9&m=10,11,12,13\\
\hline
14&m=0,1,2,3,4,5&m=6,7,8,9&m=10,11,12,13,14\\
\hline
15&m=0,1,2,3,4,5&m=6,7,8,9,10&m=11,12,13,14,15\\
\hline
16&m=0,\ds,6&m=7,8,9,10&m=11,\ds,16\\
\hline
17&m=0,\ds,7&m=8,9,10,11&m=12,\ds,17\\
\hline
18&m=0,\ds,7&m=8,9,10,11&m=12,\ds,18\\
\hline
19&m=0,\ds,8&m=9,10,11&m=12,\ds,19\\
\hline
20&m=0,\ds,9&m=10,11,12&m=13,\ds,20\\
\hline
21&m=0,\ds,10&m=11,12&m=13,\ds,21\\
\hline
22&m=0,\ds,11&m=12&m=13,\ds,22\\
\hline
23&m=0,\ds,23&&\\
\hline
24&m=0,\ds,24&&\\
\hline
25&m=0,\ds,25&&\\
\hline
26&m=0,\ds,26&&\\
\hline
27&m=0,\ds,27&&\\
\hline
28&m=0,\ds,28&&\\
\hline
29&m=0,\ds,29&&\\
\hline
30&m=0,\ds,30&&\\
\hline
\earr$$
The table above indicates that $\mu''(0)$ is always different from zero. However, we do not have a rigorous proof of this.
\end{remark}
\bibliography{fralucamax}
\bibliographystyle{abbrv}

\end{document}